\documentclass[12pt,leqno]{article}
\usepackage{amsmath,amsfonts}
\usepackage[dvips]{graphics,color}
\usepackage{latexsym}
\usepackage{amssymb}
\usepackage{graphicx}
\def\comment#1{}
\def\e0{\epsilon}

\def\f0{\phi}

\def\f{\varphi}

\def\e{\varepsilon}

\date{ }
\begin{document}
\title{The Gauss-Bonnet-Grotemeyer Theorem in spaces of constant curvature
\thanks{The project is partially supported by the grant No. 10531090 of NSFC.}}
\author{Eric L. Grinberg  and Haizhong Li }
\maketitle

\begin{abstract}
\noindent In 1963, K.P.~Grotemeyer proved an interesting variant of the Gauss-Bonnet Theorem. Let $M$ be an oriented closed surface in the Euclidean space $\mathbb R^3$ with Euler characteristic $\chi(M)$, Gauss curvature $G$ and unit normal vector field $\vec n$. Grotemeyer's identity replaces the Gauss-Bonnet integrand $G$ by the normal moment $ ( \vec a \cdot \vec n )^2G$, where $a$ is a fixed unit vector:
$
\int_M({\vec a}\cdot {\vec n})^2 Gdv=\frac{2 \pi}{3}\chi(M)
$
.
We generalize Grotemeyer's result to oriented closed even-dimesional hypersurfaces of dimension $n$ in an $(n+1)$-dimensional space form $N^{n+1}(k)$.
\end{abstract}



\section *{1. Introduction}

In 1963, K.P. Grotemeyer proved the following interesting result:

\bigskip\noindent
 {\bf Theorem 1} [Gr]) {\it Let $M$ be an oriented
closed surface in $3$-dimensional Euclidean space $\mathbb R^3$ with Gauss
curvature $G$ and a unit normal vector field $\vec{n}$. Then for any
fixed unit vector ${\vec a}$ in $R^3$, we have
$$
\int_M({\vec a}\cdot {\vec n})^2 Gdv=\frac{2 \pi}{3}\chi (M), 
\eqno (1.1)
$$
where ${\vec a}\cdot{\vec n}$ denotes the inner product of ${\vec
a}$ and ${\vec n}$, $\chi (M)$ is the Euler characteristic of
$M$.}

\medskip\noindent
{\bf Remark 1.1} Let $\{E_1,E_2,E_3\}$ be a fixed orthogonal
frame in $R^3$ and choose ${\vec a}=E_i$. We have

$$
\int_M (E_i\cdot {\vec n})^2 Gdv=\frac{2 \pi}{3}\chi(M), \quad i=1,2,3
\eqno (1.2)
$$
Noting that $\sum\limits_i (E_i\cdot {\vec n})( E_i\cdot {\vec n})={\vec n}\cdot {\vec n}=1$, we obtain the following Gauss-Bonnet formula via summation of $(1.2)$ over $i$ from $1$ to $3$:

\medskip\noindent
{\bf Corollary 1}(Gauss-Bonnet Theorem). {\it Under the same hypothesis of Theorem 1, we have}
$$
\int_M Gdv={2 \pi}\chi (M). \eqno (1.3)
$$
Thus we can consider Grotemeyer's Theorem 1 as an extended form of the Gauss-Bonnet Theorem.

\medskip
Let $n$ be even and let 
$N^{n+1}(k)$ be an $(n+1)$-dimensional simply connected Riemannian manifold of constant sectional curvature $k$. 
That is, 
$N^{n+1}(k)=\mathbb R^{n+1}$ if $k=0$; 
$N^{n+1}(k)=S^{n+1}(\frac{1}{\sqrt {k}})$, an $(n+1)$-dimensional sphere space with radius $\frac{1}{\sqrt{k}}$ if $k>0$;
 $N^{n+1}(k)=H^{n+1}(-\frac{1}{\sqrt {-k}})$, an $(n+1)$-dimensional hyperbolic space  with, as Bolyai would say, radius $\sqrt{-1}/\sqrt{k}$ if $k<0$. 
 We will often call $N^{n+1}(k)$ a {\it space form}. We will view $N^{n+1}(k)$ as standardly imbedded in an appropriate linear space 
 $L_{n+1}(k)$ ($\mathbb R^{n+2}$ if $k>0$, 
 $\mathbb R^{n+1,1}$ if $k<0$ 
 and 
 $\mathbb R^{n+1}$ if $k=0$). 
 
 This will enable us to define functions on $M$ such as 
$(\vec a\cdot \vec n )$, 
where $\vec a$ is a fixed vector in the ambient linear space, 
$\vec n$ is a normal vector field on $M$, 
and 
$( \, \cdot \, )$
denotes the inner product on the ambient linear space. The generalized Grotemeyer Theorem we have in mind can be stated as as follows:

\medskip\noindent
{\bf Theorem 2} {\it Let $n$ even, $n \ge 2$. Let $\vec x:M\to N^{n+1}(k)$ be an immersed  $n$-dimensional oriented closed hypersurface in the
$(n+1)$-dimensional space form $N^{n+1}(k)$, with Euler characteristic $\chi (M)$,Gauss-Kronecker curvature $G$ and unit normal vector field
$\vec{n}$.  Assume that $N^{n+1}(k)$ is standardly imbedded in the linear space $L_{n+1}(k)$. Then for any fixed unit vector $\vec {a}$ in
$L_{n+1}(k)$ we have
$$
\begin{array}{lcl}
{} \int_M({\vec a}\cdot {\vec n})^2 Gdv &=&
\frac{1}{n+1}[\frac{{\rm
vol}S^n(1)}{2}\chi (M)-\sum_ic_ik^i\int_MK_{n-2i}dv]
\\
{}\quad &+&\frac{k}{n+1}\int_M ({\vec a}\cdot {\vec n})({\vec a}\cdot \vec x)K_{n-1}dv-\frac{k}{n+1}\int_M ({\vec a}\cdot \vec x)^2Gdv,
\end{array}
\eqno (1.4)
$$
where the $c_i$ are constants that depend only on the dimension $n$
and $K_i$ is the $i$-th mean curvature of $M$.}

\medskip
In the case $n=2$ in the Theorem above, we obtain

\medskip\noindent
{\bf Corollary 2} {\it Let $M$ be an oriented closed surface in the
$3$-dimensional space form $N^3(k)$ with extrinsic curvature $G$
and  unit normal vector field $\vec{n}$. Then for any fixed unit vector
${\vec a}$ in the linear space $L_3(k)$  we have
$$
\begin{array}{lcl}
{}\int_M({\vec a}\cdot {\vec n})^2
Gdv &=& \frac{2 \pi}{3}\chi (M)-\frac{k}{3}{\rm vol}(M)\\
{}\quad & &+ \frac{k}{3}\int_M ({\vec a}\cdot {\vec n})({\vec a}\cdot \vec x)K_1dv -\frac{k}{3}\int_M({\vec a}\cdot \vec x)^2Gdv ,
\end{array}
\eqno (1.5)
$$
where  $K_1$ is the mean curvature of $M$ and $\chi
(M)$ is the Euler characteristic of $M$.}

\medskip\noindent
{\bf Remark 1.2} Our Corollary reduces to Grotemeyer's
original theorem in the case $k=0$.

\medskip\noindent
{\bf Remark 1.3} In the case $k=0$ and $n\geq 3$,  Theorem 2 was proved by B. -Y. Chen in [Ch] by a different method.

\bigskip\noindent
{\bf Remark 1.4}. We can recover the standard Gauss-Bonnet Theorem from our Theorem as follows. Let $m$ be the dimension of the linear space
$L_{n+1}(k)$. (Thus $m=n+1$ in the flat case, $m=n+2$ in the positive and negatively curved cases.) Let $\{E_1,\cdots, E_{m}\}$ be a fixed
orthonormal frame in $L_{n+1}(k)$; choose ${\vec a}=E_i$. Then
$$
\begin{array}{lcl}
{}
\int_M(E_i\cdot {\vec n})^2 Gdv
&=&
\frac{1}{n+1}
[\frac{{\rm vol}{S^n(1)}}{2}x(M)-\sum_i c_ik^i\int_MK_{n-2i}dv]
\\
& &+
\frac{k}{n+1}\int_M (E_i\cdot {\vec n})(E_i\cdot x)K_{n-1}dv\\
{}\quad 
& &- \frac{k}{n+1} \int_M G(E_i\cdot x)^2dv, \quad
( i=1,2,\ldots )
\end{array}
\eqno (1.6)
$$
Noting $\sum\limits_i (E_i\cdot {\vec n})( E_i\cdot {\vec n})={\vec n}\cdot {\vec n}=1$ and $\sum_i (E_i\cdot {\vec n})(E_i\cdot x)={\vec
n}\cdot x=0$, we obtain the following Gauss-Bonnet formula by summing of (1.6) over all appropriate $i$:

\bigskip\noindent
{\bf Corollary 3} (Gauss-Bonnet Theorem). {\it Under the same hypothesis of Theorem 2, we have

$$
\int_MGdv=\frac{{\rm vol}{S^n(1)}}{2}\chi(M)-\sum_i
c_ik^i\int_MK_{n-2i}dv \eqno (1.7)
$$
where $\chi(M)$ is the Euler characteristic of $M$, constants $c_i$ depends only on dimension $n$,  and $K_i$ is the $i$-th mean curvature of $M$.}

\medskip
\noindent
So we can view Theorem 2 as an extended form of the Gauss-Bonnet Theorem.

\section *{2. Reilly's operator and its properties}

In order to prove Theorem 2, we need to recall Reilly's operator and its properties.

Let $(M,g)$ be a closed $n$-dimensional Riemannian manifold, let
$\{e_1,\cdots,e_n\}$ be a local orthonormal frame field in $M$
with dual coframe field $\{\theta_1,\cdots,\theta_n\}$.
Given a symmetric tensor $\phi=\sum\limits_{i,j}\phi_{ij}\theta_i\theta_j$ defined on $M$ we
 define a second order differential operator
$$
\Box \equiv \Box_\phi : C^\infty(M)\to C^\infty (M),\qquad \Box
f=\sum\limits_{i,j}\phi_{ij}f_{ij} \eqno (2.1)
$$
where $f_{ij}$ are the components of the second covariant
differential of $f$, as follows:
$$
df=\sum\limits_i f_i\theta_i,\qquad df_i+\sum\limits_j f_j\theta_{ji}=\sum\limits_j f_{ij}\theta_j, \eqno (2.2)
$$
where $\{\theta_{ij}\}$ is the Levi-Civita connection of $g$.

\medskip\noindent
For the following criterion for self adjointness of the of the operator $\Box$ see Cheng-Yau [CY] or Li [L1],[L2].

\noindent
{\bf Proposition 2.1}  {\it Let $M$ be a closed orientable
Riemannian manifold with symmetric tensor
$\phi=\sum\limits_{i,j}\phi_{ij}\theta_i\theta_j$. Then $\Box$ is a
selfadjoint operator if and only if}
$$
\sum\limits_{j=1}^n\phi_{ij,j}=0, \quad 1\leq i\leq n. \eqno (2.3)
$$
Here $\phi_{ij,k}$ is the derivative of the tensor $\phi_{ij}$ in the direction $e_k$.

\medskip\noindent
{\bf Remark 2.1} We call $\Box$ {\it the Cheng-Yau operator}. It was
introduced by S.Y.~Cheng and S.T.~Yau in 1977 [CY]. If
$\phi=\sum\limits_{i,j}\phi_{ij}\theta_i\theta_j$ satisfies
the Cheng-Yau condition (2.3), then
$$
\Box
f=\sum\limits_{i,j}\phi_{ij}f_{ij}=\sum\limits_{i,j}(\phi_{ij}f_i)_j={\rm
div}(\phi\nabla f).
$$

Let $x:M\to N^{n+1}(k)$ be an n-dimensional closed hypersurface in
an $(n+1)$-dimensional space form of constant sectional curvature
$k$. Let$(h_{ij})$ be the components of the second fundamental
form of $M$. We recall the {\it Reilly operator}, which is a second order differential operator $L_r:
C^\infty (M)\to C^\infty (M)$ defined by
$$
L_rf=\sum\limits_{i,j}T^r_{ij}f_{ij}, \qquad f\in C^\infty (M),
\eqno (2.4)
$$
where $T^r_{ij}$ is given by
$$
T^0_{ij}=\delta_{ij},\quad T^r_{ij}=K_r\delta_{ij}-\sum\limits_k
h_{ik}T^{r-1}_{kj},\quad r=1,2,\cdots,n. \eqno (2.5)
$$
(See Reilly [Re], Rosenberg [Ro] or Barbosa-Colares [BC].)

\noindent
Denote the $r^{th}$ mean curvature of $M$ by
$$
K_r=\sum\limits_{i_1<\cdots<i_r}k_{i_1}\cdots k_{i_r},\qquad
B=(h_{ij})=(k_i\delta_{ij}). \eqno (2.6)
$$
We note that the Gauss-Kronecker curvature of $M$ is $G \equiv K_n$.

\medskip

\noindent {\bf Definition 2.1} ([Re])  The {\it r-th Newton transformation, $r\in
\{0,1,\cdots,n\}$ is the linear transformation}
$$
T_r=K_rI-K_{r-1}B+\cdots+(-1)^rB^r, \eqno (2.7)
$$
i.e.,
$$
T^r_{ij}=K_r\delta_{ij}-K_{r-1}h_{ij}+\cdots
+(-1)^r\sum\limits_{j_1,\cdots,j_r}h_{ij_1}h_{j_1j_2}\cdots
h_{j_rj}. \eqno (2.7)'
$$

If $ I \equiv i_1,\cdots, i_q$ and $J \equiv j_1,\cdots,j_q$ are multi-indices of integers between $1$
and $n$, define

$$
\delta^J_I=
\begin{cases}
 \phantom{-}1,& {\rm if}\,\,i_1,\cdots,i_q \,\text{are distinct and J is  an even} \text{ permutation of }I \\
 -1,&  {\rm if}\,\,i_1,\cdots,i_q \,\text{are distinct and J is  an odd} \text{ permutation of }I
\\
\phantom{-}0,& {\rm otherwise}.
\end{cases}
$$

%$$
%\begin{array}{lcl}
%&{}&\delta^{j_1\cdots j_q}_{i_1\cdots i_q}\\
%&=&
%\begin{cases}
%& 1, \qquad {\rm if}\,\,i_1,\cdots,i_q \,{\rm are\, distinct\,
%and\,}
%\{j_1,\cdots,j_q\} \,{\rm is\, an\, even \,permutation \,of}\, \{i_1,\cdots,i_q\}\\%
%&-1, \quad{\rm if}\,\, i_1,\cdots,i_q \,{\rm are \,distinct\, and\,}
%\{j_1,\cdots,j_q\} \,{\rm is\, an \,odd \,permutation \,of }\,\{i_1,\cdots,i_q\}\\
%&0, \qquad {\rm otherwise}.
%\end{cases}
%\end{array}
%$$
\noindent
Then we have (see Reilly [Re])
$$
K_r= \frac{1}{r!}\sum\delta^{j_1\cdots j_r}_{i_1\cdots
i_r}h_{i_1j_1}\cdots h_{i_rj_r}. \eqno (2.8)
$$

\medskip\noindent
{\bf Proposition 2.2} {\it The matrix of $T_r$ is given by}
$$
T^r_{ij}=\frac{1}{r!}\sum\delta^{j_1\cdots j_rj}_{i_1\cdots
i_ri}h_{i_1j_1}\cdots h_{i_rj_r}. \eqno (2.9)
$$

\medskip\noindent
{\bf Proposition 2.3} {\it For each $r$, we have}

(1) ${\rm div }\, T_r=\sum\limits_j T^r_{ij,j}=0$,

(2) Newton's formula:  ${\rm trace}(BT_r)=(r+1)K_{r+1}$,

\medskip

(3) ${\rm trace }(T_r)=(n-r)K_r$

\medskip\noindent
{\bf Proposition 2.4} {\it Let $\vec x:M\to N^{n+1}(k)$ be an $n$-dimensional
hypersurface with unit normal vector field ${\vec n}$. Then we have}
$$
x_i=e_i,\qquad {\vec n}_i = -\sum\limits_j h_{ij}e_j,\qquad
x_{ij}=h_{ij}{\vec n}-k x\delta_{ij}.  \eqno (2.10)
$$
$$
L_rx=(r+1)K_{r+1}{\vec n}-(n-r)kK_rx, \eqno (2.11)
$$

\medskip\noindent
{\it Proof}. Let $\vec a$ be a fixed vector in $L_n(k)$. Write
$$
f={\vec n} \cdot \vec a , \qquad g= \vec x \cdot \vec a . \eqno (2.12)
$$
Then (2.11) is equivalent to
$$
L_rg=(r+1)K_{r+1}f-(n-r)kK_rg. \eqno (2.11)'
$$

Choosing an orthonormal frame $\{e_1,\cdots,e_n, {\vec n}\}$ and
their dual frame $\{\theta_1,\cdots,\theta_n,\theta_{n+1}\}$ along
$M$ in $N^{n+1}(k)$, we have the structure equations
$$
dx=\sum_i\theta_ie_i,\quad de_i=\sum_j\theta_{ij}e_j+
\sum_jh_{ij}\theta_j{\vec n}-kx\theta_i,\quad d{\vec
n}=-\sum_{i,j}h_{ij}\theta_j e_i. \eqno (2.13)
$$
Here we have sometimes abbreviated $\vec x$ as merely $x$, for simplicity.
By use of (2.13) and through a direct calculation we get
$$
g_i= e_i \cdot \vec a ,\qquad g_{ij}=fh_{ij}-kg\delta_{ij}. \eqno (2.14)
$$

By use of proposition 2.3 and (2.14), we get
$$
L_rg=\sum_{i,j}T^r_{ij}g_{ij}=
\sum_{ij}T^r_{ij}h_{ij}f-kg\sum_{i,j}T^r_{ij}\delta_{ij}=(r+1)K_{r+1}f-k(n-r)gK_r.
$$
Thus we have proved $(2.11)'$, which is equivalent to $(2.11)$.

\noindent
Similarly, from definitions of
$f_i$, we get by use of (2.13)
$$
f_i=-\sum_j h_{ij} \, (e_j \cdot \vec a) . \eqno (2.15)
$$
Because  $\vec a$ is arbitrary, we have proved $(2.10)$ from $(2.14)$ and
$(2.15)$.

\medskip\noindent
{\bf Proposition 2.5} {\it  Let $M$ be an $n$-dimensional oriented
closed hypersurface in $(n+1)$-dimensional space form
$N^{n+1}(k)$. Then for any smooth functions $f$  and $g$ on $M$ we
have}
$$
\int_M gL_{n-1}fdv= \int_M fL_{n-1}gdv, \qquad
 \int_ML_{n-1}fdv=0. \eqno (2.16)
$$
{\it Proof}.  Choosing $r=n-1$ in (1) of proposition 2.3, and using the criteriorfrom propostion 2.1, we know that the operator $L_{n-1}$ is a  selfadjoint operator. Thus we obtain $(2.16)$.

\medskip\noindent
{\bf Proposition 2.6} {\it  Let $M$ be an $n$-dimensional
hypersurface in $(n+1)$-dimensional space form $N^{n+1}(k)$. Then
we have}
$$
G\delta_{ij}-\sum\limits_kh_{ik}T^{n-1}_{kj}=0.
\eqno (2.17)
$$

\medskip\noindent
{\it Proof}. Choosing $r=n-1$ in (2.5) and noting that $G=K_n$, we have
$$
T^n_{ij}=G\delta_{ij}-\sum\limits_k
h_{ik}T^{n-1}_{kj}.
\eqno (2.18)
$$
From the definition of $T^{n}_{ij}$ in (2.9) and the definition of
$\delta_{i_1\cdots i_ni}^{j_1\cdots j_nj}$, we have
$$
T^n_{ij}=0.
\eqno (2.19)
$$
Now (2.17) follows from (2.18) and (2.19).

\section  *{3. Proof of Theorem 2}

\medskip\noindent
{\bf Proposition 3.1}  {\it Let $x:M\to N^{n+1}(k)$ be an $n$-dimensional oriented closed hypersurface in $(n+1)$-dimensional space form
$N^{n+1}(k)$. Assume $M$ has Gauss-Kronecker curvature $G=K_n$ and a unit normal vector $\vec{n}$. Then for any fixed unit vector $\vec {a}$ in
$L_{n+1}(k)$, we have
$$
\begin{array}{rl}
0=& (n+ m) \int_M({\vec a}\cdot {\vec n})^{m+1}Gdv
-m\int_M ({\vec a}\cdot {\vec n})^{m-1}Gdv\\
 & \phantom{(n} \,\,\, -k  \int_M({\vec a}\cdot {\vec n})^m({\vec a}\cdot
\vec x)K_{n-1}dv+mk \int_M ({\vec a}\cdot {\vec n})^{m-1}({\vec a}\cdot \vec x)^2Gdv,
\end{array}
\eqno (3.1)
$$
where $K_{n-1}$ is the $(n-1)$-th mean curvature of $M$.}

\medskip\noindent
{\it Proof}. Write
$$
f=q^mx, \qquad q={\vec a}\cdot{\vec n}.
\eqno (3.2)
$$
By definition of the first derivative and the second derivative of
$f$ (see (2.2)), we have
$$
f_i=(q^m)_ix+q^mx_i,
\eqno (3.3)
$$

$$
f_{ij}=(q^m)_{ij}x+(q^m)_ix_j+(q^m)_jx_i+q^mx_{ij}.
\eqno (3.4)
$$
By definition of operator $L_{n-1}$, we have
$$
L_{n-1}(f)=x L_{n-1}(q^m)+2\sum\limits_{i,j}T^{n-1}_{ij}(q^m)_ix_j+q^m L_{n-1}x. \eqno (3.5)
$$
Let $r=n-1$ in $(2.11)$. We have
$$
L_{n-1}x=nG{\vec n}-kK_{n-1}x. \eqno (3.6)
$$
By Proposition 2.5, $(3.6)$, $(2.10)$ and proposition 2.6, we get by integrating $(3.5)$ over $M$
$$
\begin{array}{lcl}
0
&=&2\int_Mq^m(L_{n-1}x)dv+2\int_M\sum\limits_{i,j}T^{n-1}_{ij}(q^m)_ix_jdv\\
&=&2\int_M q^m(nG{\vec n}-k K_{n-1}x)dv+2\int_M \sum\limits_{i,j,k}T^{n-1}_{ij}mq^{m-1}[-h_{ik}({\vec a}\cdot e_k)]e_jdv\\
&=&
2\int_M q^m(nG{\vec n}-k K_{n-1}x)dv-2m\int_M q^{m-1}G\sum\limits_{j}({\vec a}\cdot e_j)e_jdv\\
&=& 2\int_M q^m(nG{\vec n}-k K_{n-1}x)dv-2m\int_M q^{m-1}G[{\vec
a}-({\vec a}\cdot {\vec n}){\vec n}-k({\vec a}\cdot \vec x) \vec x]dv,
\end{array}
\eqno (3.7)
$$
that is, we obtain for $m=1,2,3,\cdots$
$$
\begin{array}{cc}
   & 0=
  (n+m) \int_M q^mG{\vec n}dv -m\int_M q^{m-1}G {\vec a}dv \\
 & \phantom{(n+m}  -k \int_M
q^mK_{n-1}xdv+mk\int_M q^{m-1}(x\cdot {\vec a})Gxdv.
\end{array}
 \eqno (3.8)
$$
Taking the scalar product of $\vec a$ with both sides of (3.8), we get Proposition 3.1.

\medskip\noindent
{\bf Remark 3.1} Equation (3.1) was proved by Bang-Yen Chen in the case $k=0$ by a
 different method.

\medskip\noindent
{\bf Proof of Theorem 2} Choosing $m=1$ in Proposition 3.1, we have

$$
\begin{array}{cl}
(n+1)\int_M({\vec a}\cdot {\vec n})^2 Gdv =  &\phantom{-k}\int_M Gdv
+k\int_M ({\vec a}\cdot {\vec n})({\vec a}\cdot
\vec x)K_{n-1}dv
\\
 &-k\int_M G({\vec a}\cdot \vec x)^2dv.
\end{array}
\eqno (3.9)
$$
Because $M$ is a closed hypersurface in $N^{n+1}(k)$, the Gauss-Bonnet Theorem states in this case that
$$
\int_MGdv=\frac{{\rm vol}{S^n(1)}}{2} \chi (M)-\sum_i
c_ik^i\int_MK_{n-2i}dv,
 \eqno (3.10)
$$
where $\chi(M)$ is the Euler characteristic of $M$, the constants $c_i$
depend only on dimension $n$, and $K_i$ is the $i$-th mean
curvature of $M$. (See p. 1105 of [So]; c.f. [C1], [C2].)
Inserting (3.10) into (3.9), we have proved our Theorem 2.
On the other hand, by choosing $m=0$ in $(3.1)$, we have

\medskip

\noindent {\bf Corollary 3.1} (Bivens [Bi]) {\it Let $x:M\to N^{n+1}(k)$ be
 an $n$-dimensional closed oriented hypersurface in $N^{n+1}(k)$. Then

$$
\int_M[n({\vec a}\cdot {\vec n}) G-k({\vec a}\cdot \vec x)K_{n-1}]dv=0,
\eqno (3.11)
$$
where $\vec a$ is any fixed unit vector in the linear space $L_{n+1}(k)$, $G$ is the Gauss-Kronecker curvature of $M$ and $K_{n-1}$ is the
$(n-1)$-th mean curvature of $M$.}

\smallskip
\noindent
{\bf Remark 3.2} Write $q={\vec a}\cdot {\vec n}$, from Proposition 3.1, we have
$$
\int_M q^m Gdv\\
 =\frac{m-1}{n+m-1} [\int_M q^{m-2} Gdv - k\int_M q^{m-2}(x\cdot \vec a)^2 Gdv + \frac{k}{m-1} \int_M q^{m-1}(x \cdot \vec a ) K_{n-1} dv].
\eqno (3.12)
$$
By a direct calculation using (3.12), (3.10) and Corollary 3.1, we obtain

\medskip\noindent
{\bf Proposition 3.2} {\it Let $n$ and $m$ be even. Under the same hypothesis of Proposition 3.1, we have}
$$
\begin{array}{lcl}
&{}&
\int_Mq^mGdv\\
&=&
\frac{(m-1)(m-3)\cdots 1}{(n+m-1)(n+m-3)\cdots (n+1)}[\frac{{\rm vol}S^n(1)}{2}
\chi (M)-\sum\limits_ic_ik^i\int_MK_{n-2i}dv]\\
&{}&
-k\int_M[\frac{m-1}{n+m-1}q^{m-2}+\frac{(m-1)(m-3)}{(n+m-1)(n+m-3)}q^{m-4}+\cdots
+\frac{(m-1)(m-3)\cdots 1}{(n+m-1)(n+m-3)\cdots (n+1)}]( \vec x\cdot \vec a )^2Gdv\\
&{}&
+k\int_M[\frac{1}{n+m-1}q^{m-1}+\frac{m-1}{(n+m-1)(n+m-3)}q^{m-3}+\cdots
+\frac{(m-1)(m-3)\cdots 3}{(n+m-1)(n+m-3)\cdots (n+1)}q]( \vec x \cdot \vec a)K_{n-1}dv.
\end{array}
\eqno (3.13)
$$
{\it Also, for $n$ even and $m$ odd, we have}
$$
\begin{array}{lcl}
&{}&
\phantom{-k}\int_Mq^mGdv\\
&=& - k\int_M[\frac{m-1}{n+m-1}q^{m-2}+\frac{(m-1)(m-3)}{(n+m-1)(n+m-3)}q^{m-4}+\cdots + \frac{(m-1)(m-3)\cdots 2}{(n+m-1)(n+m-3)\cdots
(n+2)}q]( \vec x \cdot \vec a)^2Gdv
\\
&{}&
 +k\int_M[\frac{1}{n+m-1}q^{m-1}+\frac{m-1}{(n+m-1)(n+m-3)}q^{m-3}+\cdots
+\frac{(m-1)(m-3)\cdots 2}{(n+m-1)(n+m-3)\cdots n}]( \vec x \cdot \vec a )K_{n-1}dv.
\end{array}
\eqno (3.14)
$$
{\bf Note: } In the case $k=0$, Proposition 3.2 was proved by Bang-Yen Chen; see Theorem 2 in [Ch].

\bigskip\noindent
{\bf Acknowledgements}. The authors began this research works when
H. Li visited department of mathematics and statistics in
University of New Hampshire on July of 2006. H. Li would like to
thank E.L.G. and department faculty members for their hospitality
and the help they extended to him for his academic visit.

\begin{flushleft}
\medskip\noindent
\begin{tabbing}
XXXXXXXXXXXXXXXXXXXXXXXXXX*\=\kill
Eric L. Grinberg\>Haizhong Li\\
Department of Mathematics \& Statistics\>Department of Mathematical Sciences\\
University of New Hampshire  \>Tsinghua University\\
Durham, NH 03824\>100084, Beijing\\
United States of America\>People's Republic of China\\
Email: grinberg@unh.edu\>Email: hli@math.tsinghua.edu.cn

\end{tabbing}

\end{flushleft}

\vskip 0.3 in
\begin{thebibliography}{99}
\vskip 0.2in


\bibitem[BC]] J.L.M.Barbosa and A.G.Colares,  Stability of hypersurfaces with constant $r$-mean curvature,
 {\it Ann. Global Anal. Geom.}
 15(1997), 277-297.


\bibitem[Bi]] I. Bivens,  Some integral formulas for hypersurfaces in a simply connected space form,
 {\it Proc. Amer. Math. Soc.}, 88(2)(1983), 113-118.

\bibitem[Ch]] B. -Y. Chen, On an integral formula of Gauss-Bonnet-Grotemeryer, {\it Proc. Amer. Math. Soc.} 28(1971), No.1, 208-212.


\bibitem[CY]]
S.Y.Cheng and S.T.Yau,  Hypersurfaces with constant scalar
curvature, {\it Math. Ann.} 225(1977), 195-204.

\bibitem[C1]] S. S. Chern, A simple intrinsic proof of the Gauss-Bonnet formula for closed Riemannian manifolds, {\it Ann. of Math. }(2) 45(1944), 747-752.


\bibitem[C2]] S. S. Chern, On the curvatura integra in a Riemannian manifold,  {\it Ann. of Math.} (2) 46(1945), 674-684.


\bibitem[Gr]] K. P. Grotemeyer, Uber das Normalenbundel differenzierbarer Mannig faltigeiten,
{\it Ann. Acad. Sci. Fenn. Ser. A. I.} No. 336/15(1963). MR 29568.


\bibitem[L1]]
 H.Li,  Hypersurfaces with constant scalar curvature in space forms, {\it Math. Ann.} 305(1996),
665-672.

\bibitem[L2]]
H.Li,  Global rigidity theorems of hypersurface, {\it Ark. Math.}
35(1997), 327-351.

\bibitem[Re]]
R.Reilly,  Variational properties of functions of the mean
curvatures for hypersurfaces in space forms, {\it J. Diff. Geom.}
8(1973), 465-477.

\bibitem[Ro]]
H.Rosenberg,  Hypersurfaces of constant curvatures in space forms, {\it Bull. Sci. Math.} 117(1993), 211-239.


\bibitem[So]]
Gil Solanes, Integral geometry and the Gauss-Bonnet theorem in constant curvature spaces,
{\it Trans. Amer. Math. Soc.}, 358(2006), N0.3, 1105-1115.

\end{thebibliography}
\end {document}